\documentclass{article}

\usepackage{graphicx,etoolbox, xcolor
}

\usepackage{hyperref}

\usepackage{hyperref}

\usepackage{amssymb}

\newcommand{\newc}{\newcommand}

\newc{\eqnoset}{\setcounter{equation}{0}}

\newcommand{\mref}[1]{(\ref{#1})}
\newcommand{\reflemm}[1]{Lemma~\ref{#1}}
\newcommand{\refrem}[1]{Remark~\ref{#1}}
\newcommand{\reftheo}[1]{Theorem~\ref{#1}}

\newcommand{\refcoro}[1]{Corollary~\ref{#1}}

\newcommand{\refsec}[1]{Section~\ref{#1}}

\newcommand{\beq}{\begin{equation}}
	\newcommand{\eeq}{\end{equation}}
\newcommand{\beqno}[1]{\begin{equation}\label{#1}}
	
	\newcommand{\barr}{\begin{array}}
		\newcommand{\earr}{\end{array}}
	
	\newc{\bearr}{\begin{eqnarray*}}
		\newc{\eearr}{\end{eqnarray*}}
	
	\newc{\bearrno}[1]{\begin{eqnarray}\label{#1}}
		\newc{\eearrno}{\end{eqnarray}}
	
	\newc{\non}{\nonumber}
	\newc{\nol}{\nonumber\nl}
	
	\newcommand{\bdes}{\begin{description}}
		\newcommand{\edes}{\end{description}}
	\newc{\benu}{\begin{enumerate}}
		\newc{\eenu}{\end{enumerate}}
	\newc{\btab}{\begin{tabular}}
		\newc{\etab}{\end{tabular}}

	\newtheorem{theorem}{Theorem}[section]
	\newtheorem{defi}[theorem]{Definition}
	\newtheorem{lemma}[theorem]{Lemma}
	\newtheorem{rem}[theorem]{Remark}
	\newtheorem{exam}[theorem]{Example}
	\newtheorem{propo}[theorem]{Proposition}
	\newtheorem{corol}[theorem]{Corollary}
	\newtheorem{conj}[theorem]{Conjecture}

	\newcommand{\btheo}[1]{\begin{theorem}\label{#1}}
		\newc{\brem}[1]{\begin{rem}\label{#1}\em}
			\newc{\bexam}[1]{\begin{exam}\label{#1}\em}
				\newc{\bdefi}[1]{\begin{defi}\label{#1}}
					\newcommand{\blemm}[1]{\begin{lemma}\label{#1}}
						\newcommand{\bprop}[1]{\begin{propo}\label{#1}}
							\newcommand{\bcoro}[1]{\begin{corol}\label{#1}}
								\newcommand{\btheoc}[1]{\begin{conj}\label{#1}}
									\newcommand{\etheo}{\end{theorem}}
								\newc{\etheoc}{\end{conj}}
							\newcommand{\elemm}{\end{lemma}}
						\newcommand{\eprop}{\end{propo}}
					\newcommand{\ecoro}{\end{corol}}
				\newc{\erem}{\end{rem}}
			\newc{\eexam}{\end{exam}}
		\newc{\edefi}{\end{defi}}
	
	\newc{\rmk}[1]{{\bf REMARK #1: }}
	\newc{\DN}[1]{{\bf DEFINITION #1: }}
	
	\newcommand{\bproof}{{\bf Proof:~~}}
	\newc{\eproof}{{\vrule height8pt width5pt depth0pt}\vspace{3mm}}
	
	\newc{\bfrac}[2]{\dspl{\frac{#1}{#2}}}

	\newc{\nid}{\noindent}
	
	\newcommand{\dspl}{\displaystyle}
	\newc{\grad}{\nabla}
	\newc{\Div}{\mbox{div}}
	\newc{\pdt}[1]{\dspl{\frac{\partial{#1}}{\partial t}}}
	\newc{\pdn}[1]{\dspl{\frac{\partial{#1}}{\partial \nu}}}
	\newc{\pdNi}[1]{\dspl{\frac{\partial{#1}}{\partial \mathcal{N}_i}}}
	\newc{\pD}[2]{\dspl{\frac{\partial{#1}}{\partial #2}}}
	\newc{\dt}{\dspl{\frac{d}{dt}}}
	\newc{\bdry}[1]{\mbox{$\partial #1$}}
	\newc{\sgn}{\mbox{sign}}
	
	\newc{\Hess}[1]{\frac{\partial^2 #1}{\pdh z_i \pdh z_j}}
	\newc{\hess}[1]{\partial^2 #1/\pdh z_i \pdh z_j}

	\newc{\ag}{\alpha}
	\newc{\bg}{\beta}
	\newc{\cg}{\gamma}\newc{\Cg}{\Gamma}
	\newc{\dg}{\delta}\newc{\Dg}{\Delta}
	\newc{\eg}{\varepsilon}
	\newc{\zg}{\zeta}
	\newc{\thg}{\theta}
	\newc{\llg}{\lambda}\newc{\LLg}{\Lambda}
	\newc{\kg}{\kappa}
	\newc{\rg}{\rho}
	\newc{\sg}{\sigma}\newc{\Sg}{\Sigma}
	\newc{\tg}{\tau}
	\newc{\fg}{\phi}\newc{\Fg}{\Phi}
	\newc{\vfg}{\varphi}
	\newc{\og}{\omega}\newc{\Og}{\Omega}
	\newc{\pdh}{\partial}
	
	\newc{\ccG}{{\cal G}}

	\newc{\ii}[1]{\int_{#1}}
	\newc{\iidx}[2]{{\dspl\int_{#1}~#2~dx}}
	\newc{\bii}[1]{{\dspl \ii{#1} }}
	\newc{\biii}[2]{{\dspl \iii{#1}{#2} }}
	\newc{\su}[2]{\sum_{#1}^{#2}}
	\newc{\bsu}[2]{{\dspl \su{#1}{#2} }}

	\newc{\biiom}[1]{{\dspl\int_{\bdrom}~ #1 ~d\sg}}
	\newc{\io}[1]{{\dspl\int_{\Og}~ #1 ~dx}}
	\newc{\bio}[1]{{\dspl\int_{\bdrom}~ #1 ~d\sg}}
	\newc{\bsir}{\bsu{i=1}{r}}
	\newc{\bsim}{\bsu{i=1}{m}}
	
	\newc{\iibr}[2]{\iidx{\bprw{#1}}{#2}}
	\newc{\Intbr}[1]{\iibr{R}{#1}}
	\newc{\intbr}[1]{\iibr{\rg}{#1}}
	\newc{\intt}[3]{\int_{#1}^{#2}\int_\Og~#3~dxdt}
	
	\newc{\itQ}[2]{\dspl{\int\hspace{-2.5mm}\int_{#1}~#2~dz}}
	\newc{\mitQ}[2]{\dspl{\rule[1mm]{4mm}{.3mm}\hspace{-5.3mm}\int\hspace{-2.5mm}\int_{#1}~#2~dz}}
	\newc{\mitQQ}[3]{\dspl{\rule[1mm]{4mm}{.3mm}\hspace{-5.3mm}\int\hspace{-2.5mm}\int_{#1}~#2~#3}}
	
	\newc{\mitx}[2]{\dspl{\rule[1mm]{3mm}{.3mm}\hspace{-4mm}\int_{#1}~#2~dx}}
	\newc{\mitmu}[2]{\dspl{\rule[1mm]{3mm}{.3mm}\hspace{-4mm}\int_{#1}~#2~d\mu}}
	\newc{\iidmu}[2]{\iidx{#1}{#2}}
	
	\newc{\iidm}[3]{{\dspl\int_{#1}~#2~d #3}}
	
	\newc{\itQmu}[2]{\dspl{\int\hspace{-2.5mm}\int_{#1}~#2~d\mu}}
	\newc{\mitQmu}[2]{\dspl{\rule[1mm]{4mm}{.3mm}\hspace{-5.3mm}\int\hspace{-2.5mm}\int_{#1}~#2~d\mu}}

	\newc{\mitQq}[2]{\dspl{\rule[1mm]{4mm}{.3mm}\hspace{-5.3mm}\int\hspace{-2.5mm}\int_{#1}~#2~d\bar{z}}}
	\newc{\itQq}[2]{\dspl{\int\hspace{-2.5mm}\int_{#1}~#2~d\bar{z}}}

	\newc{\pder}[2]{\dspl{\frac{\partial #1}{\partial #2}}}

	\newc{\bdrom}{\bdry{\Og}}

	\newc{\bilhom}{\mbox{Bil}(\mbox{Hom}(\RR^{nm},\RR^{nm}))}
	\newc{\VV}[1]{{V(Q_{#1})}}
	
	\newc{\ccA}{{\mathcal A}}
	\newc{\ccB}{{\mathcal B}}
	\newc{\ccC}{{\mathcal C}}
	\newc{\ccD}{{\mathcal D}}
	\newc{\ccE}{{\mathcal E}}
	\newc{\ccH}{\mathcal{H}}
	\newc{\ccF}{\mathcal{F}}
	\newc{\ccI}{{\mathcal I}}
	\newc{\ccJ}{{\mathcal J}}
	\newc{\ccK}{{\mathcal K}}
	\newc{\ccP}{{\mathcal P}}
	\newc{\ccQ}{{\mathcal Q}}
	\newc{\ccR}{{\mathcal R}}
	\newc{\ccS}{{\mathcal S}}
	\newc{\ccT}{{\mathcal T}}
	\newc{\ccX}{{\mathcal X}}
	\newc{\ccY}{{\mathcal Y}}
	\newc{\ccZ}{{\mathcal Z}}
	
	\newc{\bb}[1]{{\mathbf #1}}
	
	\newc{\myprod}[1]{\langle #1 \rangle}
	\newc{\mypar}[1]{\left( #1 \right)}
	
	\newc{\BLLg}{\mathbf{\LLg}}
	
	\newc{\mA}{\mathbf{A}}
	\newc{\mB}{\mathbf{B}}
	\newc{\mC}{\mathbf{C}}
	\newc{\mD}{\mathbf{D}}
	\newc{\mE}{\mathbf{E}}
	\newc{\mF}{\mathbf{F}}
	\newc{\mJ}{\mathbf{J}}
	\newc{\mG}{\mathbf{G}}
	\newc{\mP}{\mathbf{P}}
	\newc{\mR}{\mathbf{R}}
	\newc{\mQ}{\mathbf{Q}}
	\newc{\mX}{\mathbf{X}}
	\newc{\muu}{\mathbf{u}}
	\newc{\mvv}{\mathbf{v}}
	
	\newc{\mllg}{\mathbb{\lambda}}
	\newc{\mLLg}{\mathbf{\LLg}}

	\newc{\lspn}[2]{\mbox{$\| #1\|_{\Lsp{#2}}$}}
	\newc{\Lpn}[2]{\mbox{$\| #1\|_{#2}$}}
	\newc{\Hn}[1]{\mbox{$\| #1\|_{H^1(\Og)}$}}

	\newc{\mynorm}[2]{\| #1\|_{#2}}

	\newcommand{\RR}{{\rm I\kern -1.6pt{\rm R}}}

	\newc{\itQQ}[2]{\dspl{\int_{#1}#2\,dz}}
	\newc{\mmitQQ}[2]{\dspl{\rule[1mm]{4mm}{.3mm}\hspace{-4.3mm}\int_{#1}~#2~dz}}
	\newc{\MmitQQ}[2]{\dspl{\rule[1mm]{4mm}{.3mm}\hspace{-4.3mm}\int_{#1}~#2~d\mu}}

	\newc{\MUmitQQ}[3]{\dspl{\rule[1mm]{4mm}{.3mm}\hspace{-4.3mm}\int_{#1}~#2~d#3}}
	\newc{\MUitQQ}[3]{\dspl{\int_{#1}~#2~d#3}}

	\newc{\mccP}{\mathbb{P}}
	\newc{\mccK}{\mathbb{K}}
	
	\newc{\DKTmU}{\mccK(U)}
	\newc{\DKTmUold}{(K_U(U)^{-1})^T}
	
	\newc{\myPi}{\mathbf{W}}
	\newc{\myIbar}{\bar{\ccI}_1}
	\newc{\myIhat}{\hat{\ccI}_1}
	\newc{\myIbreve}{\breve{\ccI}_0}
	
	\newc{\mmk}{\mathbf{k}}

	\newcommand{\mg}{\mathbf{g}}

	\newc{\mfu}{\mathbf{f_u}}
	\newc{\mh}{\mathbf{h}}
	\newc{\mb}{\mathbf{b}}
	\newc{\mf}{\mathbf{f}}

	\newcommand{\barrl}[2]{\barr{ll}\lefteqn{#1}\hspace{#2}&\\}

	\newc{\twomatrix}[1]{\left[\barr{cc}#1\earr\right]}
	\newc{\threematrix}[1]{\left[\barr{ccc}#1\earr\right]}

	\newc{\mN}{\mathbf{N}}
	\newc{\mI}{\mathbf{I}}
	\newc{\mH}{\mathbf{H}}

	\newc{\mk}{\mathbf{k}}
	\newc{\mr}{\mathbf{r}}

	\newc{\DIAGM}[2]{\left[\barr{ccc}#1&0\ldots&0\\
		\vdots&\ddots&\vdots\\0&\ldots0&#2\earr \right]}
	\newc{\DiagM}[2]{\mbox{diag}\left[#1
		\cdots #2 \right]}
	\newc{\vVEC}[2]{\left[\barr{c}#1\\
		\vdots\\#2\earr \right]}
	\newc{\hVEC}[2]{\left[#1
		\cdots #2 \right]}
	
	\newc{\mq}{\mathbf{q}}

	\newc{\msys}[1]{\left\{\barr{l}#1\earr
		\right.}
	\newc{\msysa}[1]{\left\{\barr{ll}#1\earr
		\right.}
	
	\newc{\bbM}{\mathbb{M}}
	\newc{\mat}[1]{\left[\barr{cc}#1\earr\right]}
	
	\newc{\me}{\mathbf{e}}

	\newc{\vecc}[2]{\left[\barr{cc}#1\\#2\earr\right]}
	\newc{\mL}{\mathbb{L}}
	
	\newc{\cO}{{\cal O}}
	
	\newc{\cM}{{\cal M}}
	
	\newc{\myega }{\eg_0(R)}
	\newc{\myeg}{\eg_1(\eg_*)}
	\newc{\myegp}{\hat{\eg}_1(\eg_*)}

	\newc{\diagA}{\mathbb{A}_d}
	\newc{\mBB}{\mathbb{B}}
	\newc{\MLT}[1]{{\cal M}_{lt}(\Og,#1)}
	\newc{\ALT}[1]{{\cal A}_{l}(\Og,#1)}
	\newc{\mM}{\mathbb{M}}
	
	\newc{\diag}[1]{\mbox{diag}(#1)}
	\newc{\off}[1]{\mbox{offdiag}(#1)}
	\newc{\mT}{\mathbb{T}}

	\usepackage{amssymb}

	\newc{\idmu}[2]{{\dspl\int_{#1}~#2~d\mu}}
	\newc{\idllg}[2]{{\dspl\int_{#1}~#2~d\llg}}
	
	\newc{\wsp}[2]{\dspl{\int}_{#1}\dspl{\int}_{#1}\frac{|#2(x)-#2(y)|^p}{|x-y|^{N+sp}}dxdy}
	\newc{\Wsp}[2]{\|#2\|_{W^{s,p}(#1)}}
	
	\newc{\dsg}{(-\Delta)^{\frac{\sg}{2}}}
	\newc{\dgeng}[1]{(-\Delta)^{\frac{#1}{2}}}
	
	\newc{\dcg}{(-\Delta)^{\frac{\cg}{2}}}
	
	\newc{\ds}{(-\Delta)^{\frac{s}{2}}}
	
	\newc{\dsab}{(-\Delta_{||})^{\frac{s}{2}}}
	\newc{\dsgab}{(-\Delta_{||})^{\frac{\sg}{2}}}

	\newc{\midx}[2]{\dspl{\rule[1mm]{3mm}{.3mm}\hspace{-4mm}\int_{#1}~#2~dx}}
	
	\newc{\midz}[3]{\dspl{\rule[1mm]{3mm}{.3mm}\hspace{-4mm}\int_{#1}~#2~ #3}}
	
	\newc{\idz}[3]{\dspl{\int_{#1}~#2~ #3}}
	
	\newc{\cfs}{C_{FS}}
	\newc{\cps}{C_{PS}}
	\newc{\nc}{\mathbf{c}}

	\newc{\xoa}{{\bf Delete?}}
	\newc{\coilai}{{\bf ???}}

\newc{\myrefc}[1]{{\bf refer #1}}
\newc{\mybook}[1]{#1}

\begin{document}

	\vspace*{-.8in}
	\begin{center} {\LARGE\em Weak local Gagliardo-Nirenberg type inequalities with a BMO term}
		
	\end{center}

	\vspace{.1in}
	
	\begin{center}

		{\sc Dung Le}{\footnote {Department of Mathematics, University of
				Texas at San
				Antonio, One UTSA Circle, San Antonio, TX 78249. {\tt Email: Dung.Le@utsa.edu}
				}}

	\end{center}

	\begin{abstract} An improvement of  a {\em (strong) global  Gagliardo-Nienberg inequality with a BMO term} is established.
		
	 \end{abstract}

\section{Introduction} \label{fracintro} In this paper we provide some improvements of  the global (strong) Gagliardo-Nienberg inequality with a BMO term in \cite{SR} and
it local (weak) Gagliardo-Nienberg inequality with a BMO term in \cite{SR,dlebook1} by extending these results in using the non-local operators fractional Laplacians. In the following, we will briefly explain the main reasons in applications that lead to these improvements.

First of all, the global (strong) Gagliardo-Nienberg inequality with a BMO term in \cite{SR}  proves that $\|Df\|_{L^{2p}(\RR^N)}^2\le C\|f\|_{BMO}\|f\|_{W^{2,p}(\RR^N)}$. This allows one to apply the interpolation inequality of Gagliardo and Nirenberg to a wide class of PDEs where one was able to have some estimate of the BMO norm of unbounded weak solutions. 
Under suitable assumptions we showed that a  weak solution $u$ to elliptic/parabolic systems will be (H\"older) regular if  the following property of $u$ is satisfied

\bdes\item[BMOsmall)] The BMO norm $\|u\|_{BMO(B_R)}$, $B_R$ denotes a ball in $\RR^N$ with radius $R$, is small if $R$ is sufficiently small.
\edes
We also elucidate in \cite{dlebook1} the pivotal role played by this property and, in some cases, the global existence of classical solutions) to the systems.

This paper has two main goals. First, we will presents an improvement of \cite{SR} by establishing a {\em Global (strong) Gagliardo-Nienberg inequality with a BMO term and fractional Laplacians.}

In \refsec{mainres}, we collect our main tools in Harmonic analysis used in our proofs. An excellent source for these is the book \cite{Grafa}.

\refsec{newWGNsec} presents the most technical proof of our main result \reftheo{GNlobalz}, a local weak  Gagliardo-Nienberg inequality. This allows us to deal with two functions ($H,u$) defined weakly in a local domain $\Og$. Several variances of this result will follow (in particular, \refcoro{locwGNcoro}).
Furthermore, we derive the global strong  Gagliardo-Nienberg inequality in \cite{SR} in \reftheo{globGNthmz} by choosing $H,u$ properly.

\section
{Notations and tools from Harmonic analysis}\label{mainres}\eqnoset

For any measurable subset $A$ of $\Og$  and any  locally integrable function $U:\Og\to\RR^m$ we denote by  $|A|$ the measure of $A$ and $U_A$ the average of $U$ over $A$. That is, $$U_A=\mitx{A}{U(x)} =\frac{1}{|A|}\iidx{A}{U(x)}.$$

We now recall some well known notions from Harmonic Analysis.

A function $f\in L^1(\Og)$ is said to be in $BMO(\Og)$ if \beqno{bmodef} [f]_{*}:=\sup_{Q\subset \Og}\mitx{Q}{|f-f_Q|}<\infty.\eeq We then define $$\|f\|_{BMO(\Og)}:=[f]_{*}+\|f\|_{L^1(\Og)}.$$

We also recall the Hardy-Littlewood maximal function
$$M(f)(x)=Mf(x):=\sup_{r>0}\midz{B_r(x)}{|f(y)|}{dy}$$
and the Hardy-Littlewood inequality if $q>1$ and $w$ belongs to the class $A_\cg$ 
\beqno{HL}\iidx{\Og}{|M(f)(x)|^qw(x)}\le C(q)\iidx{\Og}{|f(x)|^qw(x)}.\eeq

For $\cg\in(1,\infty)$ we say that a nonnegative locally integrable function $w$ belongs to the class $A_\cg$ or $w$ is an $A_\cg$ weight if the quantity
$$ [w]_{\cg} := \sup_{B\subset\Og} \left(\mitx{B}{w}\right) \left(\mitx{B}{w^{1-\cg'}}\right)^{\cg-1} \quad\mbox{is finite}.$$
Here, $\cg'=\cg/(\cg-1)$ and the supremum is taken over all cubes $B$ in $\Og$.

Here and throughout this chapter, we write $B_R(x)$ for a cube centered at $x$ with side length $R$ and sides parallel to the standard axes of $\RR^N$. We will omit $x$ in the notation $B_R(x)$ if no ambiguity can arise. We denote by $l(B)$ the side length (radius) of a cube (ball) $B$ and by $\tau B$ the cube (ball) which is concentric with $B$ and has side length (radius) $\tau l(B)$.

One of the key ingredients of our proof is the duality between the Hardy space $\ccH^1(\RR^N)$ and $BMO(\RR^N)$ space. This is the famous Fefferman-Stein theorem (see \cite{fst}). It is well known that the norm of the Hardy space can be defined by
$$\|g\|_{\ccH^1(\RR^N)}=\|g\|_{L^1(\RR^N)}+\|M_*g\|_{L^1(\RR^N)}.$$ Here, $$M_{*}g =\sup_{\eg>0}|\fg_\eg*g|$$ where  $\fg_\eg=\eg^{-N}\fg(\frac{x}{\eg})$ with $\fg\in C^\infty(\RR^N)$ which has support in the ball $B_1(y)$ centered at $y$ and $\int_{\RR^N}\fg(y)dy=1$. The definition does not depend on the choice of $\fg$ (\cite{fst}) so that throughout this paper, by the properties of $\fg$ and scaling,  we will always assume that \beqno{fgprop}0\le \fg_\eg(x)\le \eg^{-N},\; |D \fg_\eg(x)|\le \eg^{-N-1} \mbox{ for all $x\in \RR^N$}.\eeq
The last property easily comes from the diameter of $\mbox{supp}\fg_\eg$.

The famous Fefferman-Stein theorem states that there is a constant $\cfs(N)$ such that \beqno{FSineq}\left|\iidx{\RR^N}{f(x)g(x)}\right|\le \cfs(N)\|f\|_{BMO}\|g\|_{{\cal H}^1(\RR^N)}.\eeq

We will also use the definition of the  {\em centered}  Hardy-Littlewood maximal operator acting on functions $F\in L^1_{loc}(\Og)$
\beqno{maximal} M(F)(y) = \sup_\eg\{\mitx{B_\eg(y)}{F(x)}\,:\, \eg>0 \mbox{ and } B_\eg(y)\subset\Og\}.\eeq



\section{A new local weak Gagliardo-Nirenberg inequality} \label{newWGNsec} \eqnoset Let $\og$ be a function with compact support in $\RR^N$. We  assume that $D H, D u, D\og$ are  sufficiently integrable. We establish a local version of the global Gagliardo-Nirenbeg inequality in this section. The new inequality presented here is another version of  the weak GNBMO inequality in \cite{dlebook1} for fractional Laplacians ($\dsg$ with $\sg \ne1$) and also allows us to consider weak solutions and study their local regularities. The proof is also quite similar but provides some modifications and more accessible.  We also discuss several  variations of it which may be of interest in other applications.
We denote 
\newc{\mmyIbreve}{\mathbf{\myIbreve}}

\beqno{Idefz} \mI_1:=\iidx{\RR^N}{|H|^{2p}|D u|^2\og^2},\;
\mI_2:=\iidx{\RR^N}{|H|^{2p-2}|D H|^2\og^2},\eeq
\beqno{Idef2z}\mmyIbreve:=\iidx{\Og}{|H|^{2p}}.\eeq 

\newc{\Ofam}{{\mathfrak{O}_\og}}

{\bf Definition of the covering family $\Ofam$:} we cover $\Og$ by finitely many {\em sets } $B^i$'s, which has a finite intersection property, and choose $\og_i$'s such that $\og^2\le\sum_i\og^2_i$,  $\mbox{supp}\og_i\subset B^i$, and  $|D \og_i|\le \LLg_i$ for some finite numbers $\LLg_i$'s.

\btheo{GNlobalz}  
Suppose that $p\ge1$ and  $\myega, \mI_1,\mI_2,\mmyIbreve$ are finite.  For any given   $\eg_*>0$ we have a contant $C(N,p)$ such that
\beqno{GNglobalestz}\mI_1\le C(N,p)(\|u\|_{BMO(\Og_R)}+\eg_*)^2 \mI_2+E_\Ofam(u)\mmyIbreve+\iidx{\RR^N}{|H|^{2p}|D\og|^2},\eeq
where \beqno{E*def}E_\Ofam(u)=\sup_{B^i\in\Ofam}|B^i|^{-1}\iidx{B^i}{[|D u|^2+|u|^2\LLg_i^2]}.\eeq
\etheo

The estimate \mref{GNglobalestz} obviously yields the weak inequality in \cite{dlebook1} when $B^i$'s are simply balls. The quantity $E_\Ofam(u)$ depends on the geometry of $B^i$'s of $\Ofam$ which will be investigated further (to be more precise $E_\Ofam(u) \sim |B^i|^{-1}\|D (u\og_i)\|_{L^2(\RR^N)}$ ). In fact, $\og$  in our local weak inequality of \reftheo{GNlobalz} can be removed, keeping only $\mI_1, \mI_2$. However, due to the common nature of the local analysis of the regularity of weak solutions to PDE's, one has $\og$ in the consideration and \reftheo{GNlobalz} is sufficient. Furthermore, one can assume $H,u$ are smooth by using densities.

The key points are the Fefferman-Stein theorem and the Sobolev inequality.

\subsection{The proof:}

We prepare the proof by introducing the following functions, where (setting  $h=|H|^{p-1}H$ and for any $y\in \RR^N$ and $\eg>0$)
$$U_1= h D u\left(h-\mitx{B_{\eg}(y)}{h}\right),\;U_2= (h-\mitx{B_{\eg}(y)}{h})D u\mitx{B_{\eg}}{h},$$
$$U_3=D u\left(\mitx{B_{\eg}(y)}{h}\right)^2.$$

We see that $|H|^{2p}D u=h^2D u=U_1+U_2+U_3$ so that by integration by parts (writing $|D u|^2=\myprod{D u,D u}$)  \beqno{keyGNwsp}\mI_1=\iidx{\RR^N}{D(|H|^{2p}D u \og^2)u}\le \mathfrak{J}_1+\mathfrak{J}_2,\eeq where  $\mg_1,\mg_2,\mg_3$ with $\mg_i=D (U_i)$
\beqno{JJdef}\mathfrak{J}_1:=\left|\iidx{\RR^N}{(\mg_1+\mg_2) u\og^2}\right|,\; \mathfrak{J}_2:=\left|\iidx{\RR^N}{\mg_3 u\og^2}\right|+
\iidx{\RR^N}{|H|^{2p}|D u||\og||D\og|}. \eeq

Thus, $\mI_1$ can be written as $\myprod{\mathfrak{J},u}_{L^2(\RR^N)}$ which can be estimated by $\mathfrak{J}_i$'s. To estimate the  terms $\mathfrak{J}_i$'s in \mref{keyGNwsp}, we will use Fefferman-Stein theorem and  have the following lemmas.

{\em We emphasize that in the proof below we will use the Fefferman-Stein theorem to estimate $\mathfrak{J}_1$ and we will take the supremum in $\eg>0$ to do so. Meanwhile, we will not do so  for $h_{B_\eg(y)}$ in $\mg_3$ (in $\mathfrak{J}_2$). Thus, we can take any $\eg>0$ in estimating $\mathfrak{J}_2$.
}

First of all, by the Fefferman-Stein theorem again ($\mbox{supp}\og\subset\Og_R$ so that $\mbox{supp}(u\og)$ and $\mbox{supp}(\mg_i\og)\subset \Og_R\subset\Og$ so that the norms on the right hand below become local ones), $$\left|\iidx{\RR^N}{(\mg_1+\mg_2) u\og^2}\right|\le \cfs(N)\|u\|_{BMO{(\Og)}}\|(\mg_1+\mg_2)|\|_{\ccH^1(\Og_R)}.$$ 
We estimate $\|(\mg_1+\mg_2)|\|_{\ccH^1(\Og_R)}$ in the following lemmas.

\blemm{g12lemM} We have \beqno{newg1estzz}\iidx{\Og_R}{\sup_{\eg>0}|(\mg_1+\mg_2)*\fg_\eg|\og^2} \le \cps(N,\sg) \mI_1^\frac12\left\|M\left(|H|^{2(p-1)}|D H|^{2}\right)\right\|_{L^1(\Og_R)}^\frac12.\eeq
\elemm

\bproof Let us consider $\mg_1$ first and estimate the term $\iidx{\Og}{\sup_\eg|\mg_1*\fg_\eg|}$. From \mref{fgprop},  we have $|D \fg_\eg|\le \eg^{-N-1}$.
For any $y\in\Og$, we use integration by parts, ignoring $\og^2$ for simplicity  because the integrals involving $\og, D \og$ in \mref{keyGNwsp} can be treated the same way, see  \reflemm{g3lem} below. Partition again $\Og$ into finitely many {\em balls } $B^\eg$'s, which has a finite intersection property and $\mbox{supp}\fg_\eg\subset B_\eg$. Using the property of $\fg_\eg$ and then H\"older's inequality for any $s>1$, we have the following
$$\barr{lll}|\mg_1*\fg_\eg(y)|&=&\left|\iidx{B_\eg(y)}{D \fg_\eg(x-y)(h-h_{B_\eg(y)})h D u  }\right|\\&\le&\frac{C_1}{\eg}
\left|\mitx{B_\eg(y)}{|h-h_{B_\eg(y)}||hD u|}\right|
\\&\le&\frac{C_1}{\eg}
\left(\mitx{B_\eg(y)}{|h-h_{B_\eg(y)}|^s}\right)^\frac1s 
\left(\mitx{B_\eg(y)}{|hD u|^{s'}}\right)^\frac1{s'}.
\earr$$

Assuming $N>2$, we can take  $s=\frac{2N}{N-2}$ (if $N=2$, we take $s=\frac{2N}{N-1}$ in the proof) then $s_*=\frac{Ns}{N+ s}=2$. As $s_*=2$ and we can use Sobolev-Poincar\'e's inequality here. 
We have
\beqno{hestz}\frac{C_1}{\eg}
\left(\mitx{B_\eg}{|h-h_{B_\eg}|^s}\right)^\frac1s \le C\cps(N)\Psi_2^\frac{1}{2},\eeq where $\Psi_2(y)=M(|H|^{2(p-1)}|D H|^{2})(y)$. Here, we use the fact that 
$|D h| \le c(p)|H|^{p-1}|DH|, \quad h=|H|^{p-1}H.$
Setting $\Psi_3(y)=M(|hD u|^{s'})(y)$ and putting these estimates together we thus have \beqno{g1a}\sup_{\eg>0}|\mg_1*\fg_\eg| \le C
\Psi_2^\frac{1}{2}\Psi_3^\frac{1}{s'}.\eeq

By H\"older's inequality, we get
$$\iidx{\Og}{\sup_{\eg>0}|\mg_1*\fg_\eg|} \le \left(\iidx{\Og_R}{\Psi_2}\right)^\frac1{2}\left(\iidx{\Og_R}{\Psi_3^\frac{2}{s'}}\right)^\frac{1}{2}. $$

We have
$$\left(\iidx{\Og_R}{\Psi_2}\right)^\frac{1}{2} = \left\|M\left(|H|^{2(p-1)}|DH|^{2}\right)\right\|_{L^1(\Og_R)}^\frac12.$$

As $s'=\frac{2N}{N+2}<2$, we can also use \mref{HL} with the Lebesgue measure $wdx=dx$ to get from the definitions of $\mI_1, h$, $$ \left(\iidx{\Og_R}{\Psi_3^\frac{2}{s'}}\right)^\frac{1}{2}\le C\left( \iidx{\Og_R}{[|H|^{p}D u]^{2}} \right)^{1/2} = C\mI_1^{1/2} .$$

Therefore,  \beqno{g1est}\iidx{\Og_R}{\sup_\eg|\mg_1*\fg_\eg|} \le C \mI_1^\frac12\left\|M\left(|H|^{2(p-1)}|D H|^{2}\right)\right\|_{L^1(\Og_R)}^\frac12.\eeq

The term $\sup_\eg |\mg_2 *\fg_\eg|$ is handled similarly (again, replacing $h$ by its average $\mitx{B_{\eg}}{h}$ in $\Psi_3$). We then obtain an estimate like \mref{g1est} to have $$\iidx{\Og_R}{\sup_{\eg>0}|(\mg_1+\mg_2)*\fg_\eg|} \le C \mI_1^\frac12\left\|M\left(|H|^{2(p-1)}|D H|^{2}\right)\right\|_{L^1(\Og_R)}^\frac12.$$
This is \mref{newg1estzz} and the lemma is proved. \eproof

Of course, it is easy to remove the use of maximal function $\mI_2$ in \reflemm{g12lemM} and have
\blemm{g12lemr} Suppose that  $2<r$, $r'> s'$.  If we define
$$\mI_{1,r}=\iidx{\RR^N}{|H|^{(p-1)r'}(x)|\dsg u|^{r'}(x)\og^{r'}}\quad \mI_{2,r}=\iidx{\RR^N}{|H|^{(p-1)r}(x)|\dsg H|^{r}(x)\og^r}  ,$$
Then \reflemm{g12lemM}  holds in the form $$\iidx{\Og_R}{\sup_{\eg>0}|(\mg_1+\mg_2)*\fg_\eg|\og^2} \le \cps(N,\sg) \mI_{1,r}^\frac{1}{r'}\mI_{2,r}^\frac{1}{r}.$$ \elemm
\bproof Indeed, for such $r$, we can use H\"older's inequality and \mref{g1a} to have ($\Psi_2$ is defined as in the proof of \reflemm{g12lemM})
$$\iidx{\Og}{\sup_{\eg>0}|\mg_1*\fg_\eg|} \le \left( \iidx{\Og_R}{\Psi_2^\frac{r}{2}}\right)^\frac{1}{r}\left( \iidx{\Og_R}{\Psi_3^\frac{r'}{s'}}\right)^\frac{1}{r'}. $$
As $\frac{r}{2}>1$ we can use \mref{HL} for the first integral to estimated by $\mI_{2,r}$ defined here. Finally, because  $r'/s'> 1$, we apply \mref{HL} again for the second one on the right hand side to obtain the claim. \eproof

We also need to estimate $\|\mg_1+\mg_2\|_{L^1}$ but this can be done similarly. Indeed, considering the map $T:f\to (\mg_1+\mg_2)*f$ on $L^1$, we easily see that $\lim_{\eg\to0}\|T\fg_\eg\|_{L^1}\ge \|\mg_1+\mg_2\|_{L^1}$. The estimate for $\|T\fg_\eg\|_{L^1}$ can be obtained by the previous lemmas.

Similarly, to completely remove both $M,r$ from $\mI_{2,r}$ we can follow \mref{Idefz} and resume the definition of  $\mI_2$ in \mref{Idefz}
$$
\mI_2:=\iidx{\RR^N}{|H|^{2p-2}|D H|^2\og^2}.$$

\blemm{mainlocWBMOlem} There is a contants $C$ such that
\beqno{GNglobalestz00}\mI_1\le C(N,p)\|u\|_{BMO(\Og_R)}^2 \mI_2+ E_\Ofam(u)
\iidx{\Og}{|H|^{2p}}+C\iidx{\Og}{|H|^{2p}|D \og|^2}.\eeq
\elemm

\bproof First, we establish that for   $h=|H|^{p-1}H$ \beqno{GNglobalestzab}\mI_1\le C\| u\|_{BMO(\Og_R)}^2 \iidx{\RR^N}{|D h|^2\og^2}+E_\Ofam(u)
\iidx{\Og}{|h|^{2}}+C\iidx{\Og}{|H|^{2p}|D \og|^2}.\eeq

Recall that $\mg_i=D (U_i)$. For any $y\in\RR^N$ and $\eg>0$ (note that $U_i$ are the same as before)
$$U_1= hD u\left(h-\mitx{B_{\eg}(y)}{h}\right),\;U_2= (h-\mitx{B_{\eg}(y)}{h})D u\mitx{B_{\eg}}{h},$$
$$U_3=D u\left(\mitx{B_{\eg}(y)}{h}\right)^2.$$

Again as in \mref{keyGNwsp} (via integration by parts and dropping conjugates for simplicity), we have $$\barrl{\mI_1=\iidx{\RR^N}{D(|H|^{2p}Du \og^2)  u}\le }{1.5cm} &\left|\iidx{\RR^N}{(\mg_1+\mg_2)  u\og^2}\right|+\left|\iidx{\RR^N}{\mg_3  u\og^2}\right|+
\\&
\iidx{\RR^N}{|H|^{2p}|D u||\og||D\og|}. \earr$$

To estimate the last  terms, we argue similarly as before. 
First of all, by the Fefferman-Stein theorem again ($\mbox{supp}\og\subset\Og_R$), $$\left|\iidx{\RR^N}{(\mg_1+\mg_2)  u\og^2}\right|\le \cfs(N)\| u\|_{BMO{(\Og_R)}}\|(\mg_1+\mg_2)|\|_{\ccH^1(\Og_R)}.$$ 

For any $y\in\Og$, we use integration by parts, ignoring $\og^2$ for simplicity  because the integrals involving $\og, D \og$ in \mref{keyGNwsp} can be treated the same way as we estimate  $\mg_3$ in \reflemm{g3lem} later (after an integration by parts). 

Using the property of $\fg_\eg$ (the property of $\fg$ in \mref{fgprop} is changed accordingly) and then H\"older's inequality for any $s>1$, we have the following
$$\barr{lll}|\mg_1*\fg_\eg(y)|&=&\left|\iidx{B_\eg(y)}{D \fg_\eg(x-y)(h-h_{B_\eg(y)})hD u  }\right|\\&\le&\frac{C_1}{\eg}
\left|\mitx{B_\eg(y)}{|h-h_{B_\eg(y)}||hD u|}\right|
\\&\le&\frac{C_1}{\eg}
\left(\mitx{B_\eg(y)}{|h-h_{B_\eg(y)}|^s}\right)^\frac1s 
\left(\mitx{B_\eg(y)}{|hD u|^{s'}}\right)^\frac1{s'}.
\earr$$

As $s'<2$, we can use H\"older's inequality and the definition \mref{E*def} of $E_\Ofam(u)$ for the last integral (the maximal function and also its estimate \mref{HL}).Integrate the result and use Young's inequality to obtain \mref{GNglobalestz00}. The proof is complete. \eproof

Concerning the last term $\mathfrak{J}_2$ involving $\mg_3$ in \mref{keyGNwsp}, we have the following lemma. Notice that we do not have to estimate ${\cal H}^1$ norms although there is $\eg$ in the definition of $\mg_3$.

\blemm{g3lem} For any $\eg, \thg>0$ we have
$$\mathfrak{J}_2\le E_\Ofam(u)\iidx{\Og}{h^2}+\thg\iidx{\Og}{|H|^{2p}|D u|^2\og^2}+C(\thg)\iidx{\Og}{|H|^{2p}|D \og|^2}.$$
where (see \mref{E*def} and the definition of the family $\Ofam$)  $$E_\Ofam(u)=\sup_{B^i\in\Ofam}|B^i|^{-1}\iidx{B^i}{[|D u|^2+|u|^2\LLg_i^2]}.$$
\elemm

\bproof From  the definition of $\Ofam$, we first have by integration by parts on each $B^i$ and then Young's inequality that
$$\barr{lll}\left|\iidx{\Og}{\mg_3 u\og^2}\right|&\le& \sum_i\left| \iidx{B^i}{|h_{B^i}|^2D u u\og_i^2}\right|\\
&\le& \sum_i\iidx{B^i}{|h_{B^i}|^2[|D u|^2\og_i^2+|u|^2|\LLg_i|^2]}.\earr$$

The right hand side is estimated by $$\sum_i|h_{B^i}|^2\iidx{B^i}{[|D u|^2\og_i^2+|u|^2\LLg_i^2]}\le \sup_i\iidx{B^i}{[|D u|^2\og_i^2+|u|^2\LLg_i^2]}\sum_i|h_{B^i}|^2.$$
Using the inequality $|h_{B^i}|^2\le (h^2)_{B^i}$ and the finite intersection property of $B^i$'s we bound the last term by
\beqno{g3esta} \sup_i|B^i|^{-1}\iidx{B^i}{[|\dsg u|^2\og_i^2+|u|^2\LLg_i^2]}\iidx{\Og}{h^2}.\eeq

Regarding the second term is $\mathfrak{J}_2$, we use Young's inequality
$$\iidx{\Og}{|H|^{2p}|D u||\og||\dsg\og|}\le \iidx{\Og}{|H|^{2p}[\thg|D u|^2\og^2+C(\thg)|D \og|^2]}.$$

Combining the above estimates, we obtain the lemma. \eproof

\brem{J2rem} The reason we did not use a partition in estimating the second term in $\mathfrak{J}_2$ is that we will obtain a better estimate because the number of the balls $B^i$'s would be very large. Thus, leaving $\og$ as it is would be a better choice. Also, we separate $E_\Ofam$ from the second term in $\mathfrak{J}_2$ to freely choose $\Ofam$ from $\og$ later on. It is important to note that the factor $E_\Ofam$ is multiplied by the  integral of $|h|^2=|H|^{2p}$ over $\Og$. 
\erem

Finally, we are ready to have

{\bf Proof of \reftheo{GNlobalz}:}
The last term in \mref{keyGNwsp} is treated  in \reflemm{g3lem}. Putting the estimates in \reflemm{mainlocWBMOlem} and \reflemm{g3lem} together , we obtain 
\reftheo{GNlobalz}. \eproof

\brem{constantC(N,p)} Inspecting the proof of \reflemm{g12lemM} and \reflemm{g12lemr} and the definitions of $h,\Psi_2, \Psi_3$, we remark here that $C(N,p) \sim C(N)p$ as we use the product rule for $h=H^p$ once  to estimate $\mathfrak{J}_1$'s. Here, the constant $C(N)\sim (C_{PS}(N)\cfs(N))^2$ (where $C_{PS}, \cfs$ are the constants in the Poincar\'e-Sobolev and Fefferman-Stein (see \mref{FSineq}) inequalities). \erem

\brem{EOfam}As the choices of  $\thg$ and $\Ofam$ are independent, we will later choose $B_i$ large so that $\thg, |B_i|^{-1}, C(\thg)|D\og_i|\to 0$  to obtain a better version of this theorem (note that the supports of $H,u$ can still be bounded).
We eventually prove that $E_{\Ofam}\to0$  and obtain a much better version of \reftheo{GNlobalz} in the corollaries below.\erem

Recall that (see \mref{E*def} and the definition of the family $\Ofam$)  $$E_\Ofam(u)=\sup_{B^i\in\Ofam}|B^i|^{-1}\iidx{B^i}{[|D u|^2+|u|^2\LLg_i^2]}.$$

We now look carefully at the following two terms in $\mathfrak{J}_2$ and the choices of $\Ofam, \og$ $$E_\Ofam(u)\iidx{\Og}{|H|^{2p}} \mbox{ and } \iidx{\Og}{|H|^{2p}|D \og|^2}.$$

Assume that $u\in BMO(B_R)$, $k<1$, and we define $\og$ to be a cut-off function for $B_R, B_{2R}$. Of course when $R\to\infty$, we see that $|D\og|, \LLg, |B_R|^{-1}\to 0$. Therefore if we choose $\Ofam$ to be just the ball $B_R$ then we have $E_\Ofam(u)\to 0$. We can let $\thg\to0$ too.
We thus have a local and global weak Gagliardo-Nirenberg type result as follows

\bcoro{locglobWGNcoro} Assume that $R>0$ and measurable $H,u$ on $B_R$ are such that $u\in BMO(B_R)$ and $H, u$ can be extended to $\RR^N$ such that $u, D u\in L^2(\RR^N),\; |H|^{2p} \in L^1(\RR^N)$. Then,
\beqno{GNglobalestzcoro}\iidx{\RR^N}{|H|^{2p}|D u|^2}\le C(N,p)\|u\|_{BMO(B_R)}^2\iidx{\RR^N}{|H|^{2p-2}|D H|^2}.\eeq
The above integrals can be both infinite. If $|H|^{2p-2}|D H|^2 \in L^1(\RR^N)$ then they are both finite.
\ecoro

We can not replace the integral on the right hand side  over $\RR^N$ by the one over $B_R$ (even $H\in W^{1,2}(B_R)\cap L^{2p-2}(B_R)$). Of course, the left hand side integral can be reduced to that over $B_R$. A reason is that if  $H\in W^{1,2}(B_R)\cap L^{2p-2}(B_R)$ then we can extend $H$ to $\RR^N$ and  have a constant $C$ such that
$$\iidx{\RR^N}{|H|^{2p-2}}\le C\iidx{B_R}{|H|^{2p-2}},\; \iidx{\RR^N}{|D H|^{2}}\le C\iidx{B_R}{|D H|^{2}}$$ but we can not find such a constant $C$ and
$$\iidx{\RR^N}{|H|^{2p-2}|D H|^2}\not\le C\iidx{B_R}{|H|^{2p-2}|D H|^2}.$$

If $u\in BMO(B_R)$ and $u, |H|^{p}H\in W^{1,2}_0(B_R)$ then by extending $H$ to be zero outside $B_R$ we obtain
$$\iidx{B_R}{|H|^{2p}|D u|^2}\le C\|u\|_{BMO(B_R)}^2\iidx{B_R}{|H|^{2p-2}|D H|^2}.$$
In fact, this is also true if $|H|^{p-1}H\in W^{1,2}(B_R)$ as  the last integral is in fact that of $|D (|H|^{p-1}H)|^2$ if we inspect the proof of  \reflemm{g12lemM} where we estimate $\mg_1$. Thus, $\og$ and $\mathfrak{J}_2$ in our local weak inequality of \reftheo{GNlobalz} can be removed, keeping only $\mI_1, \mI_2$, if $\Og$ is an extended domain (as $B_R$). However, due to the common nature of the local analysis of the regularity of weak solutions to PDE's, one has $\og$ in the consideration and \reftheo{GNlobalz} is sufficient. 

The fact that $E_\Ofam\to0$, when $\Ofam=\{B_R\}$ with $R\to\infty$, greatly improves \reftheo{GNlobalz} so that we state it as follows (keep in mind that $H,u$ are representatives of classes).

\bcoro{locwGNcoro}  If $u\in BMO(\Og)$, $u \in W^{1,2}(\Og)$, $|H|^{p-1}H\in W^{1,2}(\Og)$, and $\Og$ is an extended domain, then $$\iidx{\Og}{|H|^{2p}|D u|^2}\le C(N,p)\|u\|_{BMO(\Og)}^2\||H|^{p-1}D H\|_{L^{2}(\Og)}^2.$$
The above inequality also holds when $\Og=\RR^N$. 
\ecoro

\brem{constantC(N,p)a} For future references, by \refrem{constantC(N,p)}, note that $C(N,p) \sim C(N)p$, the constant $C(N)\sim (C_{PS}(N)\cfs(N))^2$ (where $C_{PS}, \cfs$ are the constants in the Poincar\'e-Sobolev and Fefferman-Stein inequalities). \erem

If $|H|^{p-1}H, |H|^{p-1}D H\in L^2(\Og)$, then $|H|^{p-1}H\in W^{1,2}(\Og)$. The above corollary also  yields
$$\iidx{\Og}{|H|^{2p}|D u|^2}\le C(N,p,\Og) \|u\|_{BMO(\Og)}^2\iidx{\Og}{|H|^{2p-2}|D H|^2}.$$

The integral on the left hand side is finite. This assertion is not trivial under the assumptions that $u\in BMO(\Og)$ and $ |H|^{p-1}H\in W^{1,2}(\Og)$ (even when $H=u$). When $p=1$, we also have
$$\mbox{$u\in BMO(\Og)$, and $ H\in W^{1,2}(\Og)$}\Rightarrow\iidx{\Og}{|H|^{2}|D u|^2}\le C(N,\Og)\|u\|_{BMO(\Og)}^2\|H\|_{W^{1,2}(\Og)}^2.$$

This leads to the special case when $H=u$. We can easily iterate \refcoro{locwGNcoro} to assert that:
$$\mbox{If  $H\in W^{1,2}(\Og)\cap BMO(\Og)$,  then for all $p\ge1$ }\iidx{\Og}{|H|^{2p}|D H|^2}<\infty.$$

As a byproduct of \refcoro{locwGNcoro}, we also have the following theoretic result. If $\Og$ is a bounded extension domain and $u\in BMO(\Og)$ and $H\in W^{1,2}(\Og)$, then by H\"older's inequality for $p>1$
$$\iidx{\Og}{|u|^2|D H|^2}\le \||u|^2\|_{L^{p'}(\Og)} \||D H|^2\|_{L^{p}(\Og)}.$$

By the continuity of the map $p\to \|g\|_{L^p(\Og)}$ for a given $g\in L^p(\Og)$ and if $u\in BMO(\Og)$ then $u\in  L^q(\Og)$ for all $q>1$, we conclude that
$$\iidx{\Og}{|u|^2|D H|^2}\le C(\|u\|_{BMO(\Og)}, \|D H\|_{L^{2}(\Og)})<\infty.$$
Together with \refcoro{locwGNcoro}, we now see that $uH\in W^{1,2}(\Og)$. Thus,

\bcoro{Walgeb} Assume that $\Og$ is a bounded extension domain and that $u,H\in W^{1,2}(\Og)$. If either $u\in BMO(\Og)$ or $H\in BMO(\Og)$, then $uH\in W^{1,2}(\Og)$. 
\ecoro

By iterating \mref{GNglobalestzcoro}  in $p$, if $u\in BMO$ and $|H|^p, D H\in L^2(\RR^N)$, then $|H|^pD u\in L^2(\RR^N)$ for all $p\ge1$. Hence, (with $H=u$) we have for $H\in BMO(\RR^N)\cap W^{1,2}(\RR^N)\cap L^2(\RR^N)$
\beqno{GNglobalestzcoroHisu}\iidx{\RR^N}{|H|^{2p}|D H|^2}\le C(N,p)\|H\|_{BMO(\RR^N)}^2 \iidx{\RR^N}{|H|^{2p-2}|D H|^2}. \eeq
This trivial if $H\in L^\infty(\RR^N)\cap W^{1,2}(\RR^N)$.

We also have from \mref{GNglobalestzcoro} that
\bcoro{GNlobalzcoro1}   Assume $p\ge1$,  and  $u\in BMO(\RR^N)$. If $|H|^{p-1}D H\in L^2(\RR^N)$, then there is a contants $C$ such that
\beqno{GNglobalestzcoro1}\iidx{\RR^N}{|H|^{2p}|D u|^2}\le C(N,p)\|u\|_{BMO(\RR^N)}^2 \iidx{\RR^N}{|H|^{2p-2}|D H|^2}. \eeq
\ecoro
When $H=u$ we have for all $p\ge1$ that
\beqno{GNglobalestzcoro2}\iidx{\RR^N}{ |H|^{2p}|D H|^2}\le C(N,p)\|H\|_{BMO(\RR^N)}^2 \iidx{\RR^N}{ |H|^{2p-2}|DH|^2}. \eeq

\subsection{$H=D u$ and \cite[Theorem 1.1]{SR}} 
The following result is just \refcoro{locwGNcoro} when $\Og=\RR^N$ and $H=D u$ .
\btheo{globGNthmz} For $p\ge1$, there is a constant $C(N,p)$ such that
\beqno{Idefzweakglob} \iidx{\RR^N}{|D u|^{2p+2}} \le C(N,p)\|u\|_{BMO}^2\iidx{\RR^N}{|D u|^{2p-2}|D^2u|^2}.\eeq This also holds when $\RR^N$ is replaced by an extended domain $\Og$ and $C(N,p)=C(N,p,\Og)$.
\etheo

By H\"older's inequality, we have
$$\iidx{\RR^N}{|D u|^{2p-2}|D^2u|^2}\le \left(\iidx{\RR^N}{|D u|^{2p+2}}\right)^\frac{p-1}{p+1}\left(\iidx{\RR^N}{|D^2 u|^{p+1}}\right)^\frac{2}{p+1}.$$

Using this in \mref{Idefzweakglob} and canceling, we get
$$\iidx{\RR^N}{|D u|^{2p+2}} \le C(N,p)\|u\|_{BMO}^2\iidx{\RR^N}{|D^2u|^{p+1}}.$$
This is exactly \cite[Theorem 1.1]{SR} when $p\ge2$ ($p\to p+1$).

\bibliographystyle{plain}

\begin{thebibliography}{10}
	
	
\bibitem{fst}  C. Fefferman and E. M. Stein, $H^p$ spaces of several variables, {\em Acta Math.} 129 (1972)
137–193

	
	\bibitem{Gius}
	E. Giusti.
	\newblock {\em Direct Methods in the Calculus of Variations}.
	\newblock World Scientific, Singapore, 2003.
	
	\bibitem{Grafa} L. Grafakos. \newblock {\em Modern Fourier Analysis (Graduate Texts in Mathematics)}. Springer; 2nd ed. 2009 edition (November 26, 2008).
	
	
	\bibitem{dlebook} D. Le, \newblock{\em Strongly Coupled Parabolic and Elliptic Systems: Existence and Regularity of Strong/Weak Solutions.} De Gruyter, 2018.
	
	\bibitem{dlebook1} D. Le, \newblock{\em Cross Diffusion Systems: Dynamics, Coexistence and Persistence.} De Gruyter, 2022.
	
	
	
	
	
	
	\bibitem{SR} P. Strzelecki.  \newblock{Gagliardo Nirenberg inequalities with a BMO term.} {\em Bull. London Math. Soc.} Vol. 38, pp. 294-300, 2006.
	








\end{thebibliography}

\end{document}